\documentclass{article}
\usepackage{mathrsfs}
\usepackage{amsmath}
\usepackage{amssymb, amsmath}
\usepackage{graphicx}
\hoffset - 10 mm

\catcode`\@=11 \@addtoreset{equation}{section}

\catcode`\@=12
\usepackage{color}

\newtheorem{Theorem}{Theorem}[section]
\newtheorem{Proposition}{Proposition}[section]
\newtheorem{Lemma}{Lemma}[section]
\newtheorem{Corollary}{Corollary}[section]
\newtheorem{Definition}{Definition}[section]

\newtheorem{Remark}{Remark}[section]

\newcommand{\newcom}{\newcommand}
\newcommand{\bTheorem}[1]{
\begin{Theorem} \label{T#1} }
\newcommand{\eT}{\end{Theorem}}

\newcommand{\bProposition}[1]{
\begin{Proposition} \label{P#1}}
\newcommand{\eP}{\end{Proposition}}

\newcommand{\bLemma}[1]{
\begin{Lemma} \label{L#1} }
\newcommand{\eL}{\end{Lemma}}

\newcommand{\bCorollary}[1]{
\begin{Corollary} \label{C#1} }
\newcommand{\eC}{\end{Corollary}}

\newcommand{\beq}{\begin{equation}}
\newcommand{\eeq}{\end{equation}}
\newcom{\ben}{\begin{eqnarray}}
\newcom{\een}{\end{eqnarray}}
\newcom{\beno}{\begin{eqnarray*}}
\newcom{\eeno}{\end{eqnarray*}}
\newcom{\bali}{\begin{aligned}}
\newcom{\eali}{\end{aligned}}

\newcommand{\bFormula}[1]{
\begin{equation} \label{#1}}
\newcommand{\eF}{\end{equation}}

\newcommand{\f}{\frac}
\newcommand{\Om}{\Omega}

\newcommand{\p}{\partial}

\newcommand{\vr}{\varrho}

\newcommand{\vt}{\vartheta}
\newcommand{\vu}{\vc{u}}

\newcommand{\vv}{\vc{v}}

\newcommand{\vc}[1]{{\boldsymbol #1}}
\newcommand{\Div}{{\rm div}}
\newcommand{\Grad}{\nabla}

\newcommand{\dx}{{\rm d} x}

\newcommand{\dt}{{\rm d} t }
\newcommand{\ds}{{\rm d} s}

\newcommand{\dxdt}{\dx\dt}

\newcommand{\ep}{\varepsilon}

\font\F=msbm10 scaled 1000
\newcommand{\R}{\mbox{\F R}}

\newcommand\Cbox[2]{%
    \newbox\contentbox%
    \newbox\bkgdbox%
    \setbox\contentbox\hbox to \hsize{%
        \vtop{
            \kern\columnsep
            \hbox to \hsize{%
                \kern\columnsep%
                \advance\hsize by -2\columnsep%
                \setlength{\textwidth}{\hsize}%
                \vbox{
                    \parskip=\baselineskip
                    \parindent=0bp
                    #2
                }%
                \kern\columnsep%
            }%
            \kern\columnsep%
        }%
    }%
    \setbox\bkgdbox\vbox{
        \color{#1}
        \hrule width  \wd\contentbox %
               height \ht\contentbox %
               depth  \dp\contentbox
        \color{black}
    }%
    \wd\bkgdbox=0bp%
    \vbox{\hbox to \hsize{\box\bkgdbox\box\contentbox}}%
    \vskip\baselineskip%
}

\begin{document}


\title{\bf On global-in-time weak solutions to the magnetohydrodynamic system of compressible inviscid fluids}

\author{Eduard Feireisl \\ Faculty of Mathematics and Physics, \\ Charles University, Sokolovsk\'a 83, CZ-186 75 Prague 8, Czech Republic
\\ and \\ TU Berlin, Strasse des 17. Juni, Berlin, Germany\\ feireisl@math.cas.cz\\ \\
Yang Li \\ Department of Mathematics, \\ Nanjing University, Hankou Road 22, 210093, Nanjing, China \\ lymath@smail.nju.edu.cn\\
}

\maketitle

{\centerline {\bf Abstract }}
{We consider the motion of an inviscid compressible fluid under the mutual interactions with magnetic field.
We show that the initial value problem  is ill--posed in the class of weak solutions for a large class of physically
admissible data. We also consider the same problem for inviscid heat--conductive fluid and show the same result under
certain restrictions imposed on the magnetic field. The main tool is the method of convex integration adapted to the
Euler system with ``variable coefficients''.}

{\bf Keywords: }{magnetohydrodynamic system, compressible flow, weak solutions, convex integration}

{\bf Mathematics Subject Classification.} 76W05, 76N15, 35D30.

\section{Introduction}
The time evolution of electrically conducting inviscid compressible fluid interacting with a magnetic field is described by the system of magnetohydrodynamics (MHD). The conservation of mass, the balance of momentum and the Maxwell system for the magnetic
field read as (see \cite{HC}):
\begin{equation}\label{lp1}
\left\{\begin{aligned}
& \p_t\vr + \Div_x(\vr\vv)=0,\\
& \p_t(\vr\vv)+\Div_x(\vr\vv \otimes \vv)+\Grad_x{p(\vr)}=
\mathbf{curl}_x \mathbf{B} \times \mathbf{B},\\
& \p_t \mathbf{B}=\mathbf{curl}_x (\vv \times \mathbf{B})-\nu \mathbf{curl}_x \mathbf{curl}_x \mathbf{B},\\
& \Div_x \mathbf{B}=0.\\
\end{aligned}\right.
\end{equation}
The unknowns are the fluid density $\vr = \vr(t,x)$, the velocity field $\vv = \vv(t,x) \in \R^3$, and the magnetic field $\mathbf{B}
= \mathbf{B}(t,x) \in \R^3$, depending on the time variable $t\in \R ^+$ and the space variable $x=(x_1,x_2,x_3)\in \R^3$.
The symbol $p(\vr)$ stands for the pressure, $\nu>0$ is the resistivity coefficient which acts as the magnetic diffusion. We normalize $\nu=1$ in what follows.

Noticing that
\[
\mathbf{curl}_x \mathbf{B} \times \mathbf{B}=\Div_x (\mathbf{B}\otimes \mathbf{B})-\Grad_x{\left(\f{1}{2}|\mathbf{B}|^2\right)},
\]
\[
\mathbf{curl}_x \mathbf{curl}_x \mathbf{B}=\Grad_x \Div_x \mathbf{B}- \Delta \mathbf{B},
\]
we may rewrite the system (\ref{lp1}) in the form
\begin{equation}\label{lp2}
\left\{\begin{aligned}
& \p_t\vr + \Div_x(\vr\vv)=0,\\
& \p_t(\vr\vv)+\Div_x(\vr\vv \otimes \vv-\mathbf{B}\otimes \mathbf{B})+\Grad_x\left(p(\vr)+\f{1}{2}|\mathbf{B}|^2\right)=
\mathbf{0},\\
& \p_t \mathbf{B}=\mathbf{curl}_x (\vv \times \mathbf{B})+\Delta \mathbf{B},\\
& \Div_x \mathbf{B}=0.\\
\end{aligned}\right.
\end{equation}

There is a well--developed mathematical theory of \emph{incompressible} MHD equations (i.e., $\Div_x\vv=0$). In the case of viscous magnetically resistive fluids, the existence, uniqueness and large time behavior of strong and weak solutions have been studied by Duvaut and Lions \cite{DL}, Sermange and Temam \cite{ST}, among others. We refer to Cao and Wu \cite{CW} for global regularity of MHD equations with mixed partial dissipation and magnetic diffusion in two space dimensions. The problem is much more involved when there is only viscosity or only magnetic diffusion present. In the case of viscous fluids without magnetic diffusion, global-in-time existence of small classical solutions has been obtained by Lin et al. \cite{LXZ} in two space dimensions, Xu and Zhang \cite{XZ} in three space dimensions. See also \cite{TW,PZZ} for related results on the initial-boundary value problem. The problem with zero viscosity and
with magnetic diffusion was studied by Cao and Wu \cite{CW} in the 2-D case, where the existence of global-in-time weak solution with $H^1$ initial data was established. See also \cite{ZZ} for the existence of global-in-time classical solution by requiring smallness and certain symmetries of the initial data. As for the inviscid and non-resistive case, we refer to Bardos et al. \cite{BSS} for the existence and large time behavior of global-in-time classical solution. By adapting the arguments of convex integration of De Lellis and Sz\'ekelyhidi \cite{DS1}, Bronzi et al. \cite{BLN} proved the existence of global-in-time weak solution to the symmetry reduced MHD equations with compact support in space and time. We remark that proving the non-uniqueness of weak solutions in the context of \emph{incompressible} flows stems from Scheffer \cite{Sche} and Shnirelman \cite{Shir} on Euler equations. See also \cite{CM1,CFG,DS2,Shvy} and the references therein for more results on non-uniqueness of weak solutions in the context of \emph{incompressible} flows. We refer to \cite{CKS} for the conservation of energy and magnetic helicity to ideal MHD equations under suitable assumptions imposed on weak solutions (see also \cite{FL}).

The theory of \emph{compressible} MHD fluid flows is more involved. The case of a viscous magnetically resistive fluid
was studied in \cite{DF}, where the existence of global--in-time weak solution to the full MHD system
is obtained (see also Hu and Wang \cite{HW} for the isentropic regime). In \cite{LXZ2,PG}, results on the existence and uniqueness of global-in-time classical solution have been obtained under smallness assumptions imposed on the initial data.
In the case of viscous non-resistive fluids, we refer to Wu and Wu \cite{WW}, Tan and Wang \cite{TW} for global-in-time existence and uniqueness of small classical solution; see also \cite{LIS} for the existence of global-in-time weak solutions with certain symmetry in two space dimensions. Fan et al. \cite{FJN} established the existence of global-in-time weak solution to the resistive planar MHD system with zero shear viscosity. Recently, Gwiazda et al. \cite{GMS} obtained a sufficient condition for the energy conservation for weak solutions to inviscid non-resistive \emph{compressible} MHD equations.

The problem of
global--in--time solvability of the MHD system for an inviscid, compressible and resistive or non--resistive fluid in the space dimension
$N=2,3$ remains largely open. Klingenberg and Markfelder \cite{KM} used the piece--wise constant data ansatz,
similar to Luo, Xie and Xin \cite{LXX}, to show ill--posedness for
the inviscid non--resistive model under simplified symmetry hypothesis by using the method of convex integration. In this context,
it is worth mentioning a very interesting recent result of Dai \cite{Dai}, where ill posedness is established for the
3-D viscous MHD equations with Hall effect.

In this paper, we consider the problem of global well/ill posedness
of the MHD system, where the fluid is compressible, inviscid but still magnetically resistive. To simplify presentation, we impose the space periodic boundary conditions and restrict ourselves to the physically relevant 3-D geometry, meaning the physical space is the (flat) three-dimensional torus
\begin{equation} \label{lp4b}
\Om:= \left( [0,1]|_{\{0,1\}} \right)^3.
\end{equation}
Problem \eqref{lp2}, \eqref{lp4b} is supplemented by the initial conditions:
\beq\label{lp4}
\vr(0,\cdot)=\vr_0,\,\, \vv(0,\cdot)=\vv_0,\,\,\mathbf{B}(0,\cdot)=\mathbf{B}_0.
\eeq

Our goal is to show that the problem \eqref{lp2}--\eqref{lp4} is globally solvable but essentially ill posed in the class of weak solutions. Specifically we establish the following results:
\begin{itemize}

\item Given sufficiently regular initial data, the MHD system \eqref{lp2}--\eqref{lp4} admits infinitely many global--in--time weak solutions,
see Theorem \ref{eom1}.

\item There is a vast class of initial data for which the MHD system \eqref{lp2}--\eqref{lp4}
admits infinitely many global--in--time weak solutions satisfying the relevant energy balance, see Theorem \ref{eom2}.

\item We extend the previous results to the case of heat conducting fluid under certain symmetry restrictions, see Theorem \ref{eom4}.

\end{itemize}

The paper is organized as follows.
In Section \ref{P}, we introduce the concept of weak solution to the MHD system and formulate our basic ill posedness result
in the general setting.
In Section \ref{ref}, we reformulate the MHD system in a convenient form and introduce the definition of subsolutions. The proof of
the basic ill posedness result is then finished by means of a suitable version of oscillatory lemma and the convex integration scheme
in Section \ref{osc}. In Section \ref{conrm}, we establish the existence of infinitely many weak solutions to (\ref{lp2}) that satisfy a relevant form of
energy inequality. Extension to the case of a heat-conductive fluid is presented under a special symmetry assumption, see Section \ref{exth}.

\section{Preliminaries, the main result in the general case}
\label{P}

We start by introducing the concept of weak solution to the MHD system (\ref{lp2})-(\ref{lp4}).

\begin{Definition}\label{defn1}
Let $\Omega$ be the flat torus introduced in \eqref{lp4b}.
A trio $[\vr,\vv,\mathbf{B}]$ is said to be a \emph{weak solution} to \eqref{lp2}--\eqref{lp4} in the time-space domain $(0,T)\times \Om$
if:
\begin{itemize}
\item {$\vr(t,x)\geq 0$  for a.e. $(t,x)\in (0,T)\times \Om$;}

\item { $\int_0^T \int_{\Om}\left(\vr \p_t \phi +\vr \vv\cdot \Grad_x \phi\right) \dxdt +\int_{\Om}\vr_0 \phi(0,\cdot)\dx=0$
      for any $\phi\in C_c^{\infty}( [0,T)\times \Om)$;}

\item {  $\int_0^T \int_{\Om}\left(\vr \vv\cdot\p_t \vc{\phi} +\left(\vr \vv\otimes \vv-\mathbf{B}\otimes \mathbf{B}\right): \Grad_x \vc{\phi}+\left(p(\vr)+\f{1}{2}|\mathbf{B}|^2\right)\Div_x \vc{\phi}\right) \dxdt $\\
    $+\int_{\Om}\vr_0 \vv_0\cdot \vc{\phi}(0,\cdot)\dx=0$
      for any $\vc{\phi}\in C_c^{\infty}( [0,T)\times \Om;\R^3)$;}

\item {$\int_0^T \int_{\Om}\left[ \mathbf{B}\cdot \p_t \vc{\phi}+(\vv \times \mathbf{B})\cdot  \mathbf{curl}_x \vc{\phi} -\Grad_x \mathbf{B}: \Grad_x \vc{\phi}     \right] \dxdt +\int_{\Om}\mathbf{B}_0 \cdot\vc{\phi}(0,\cdot)\dx=0$\\
      for any $ \vc{\phi}\in C_c^{\infty}( [0,T)\times \Om;\R^3)$;}

\item{$ \int_{\Om}\mathbf{B}\cdot \Grad_x \phi \dx =0$
    for any $\phi\in C^{\infty}( \Om)$, a.e. $t\in (0,T)$.
    }
\end{itemize}
\end{Definition}

Our first result concerns existence and non-uniqueness of global-in-time weak solutions to the problem (\ref{lp2})--(\ref{lp4}) for any smooth initial data.

\begin{Theorem}\label{eom1}
Let $T\in (0,\infty)$ be given. Assume that
\[
p(\cdot)\in C([0,\infty))\cap C^\infty((0,\infty)),\,\, p(0)=0,
\]
\[
\vr_0\in C^3(\Om),\,\,\vr_0>\underline{\vr}>0,\,\,\vv_0 \in C^3(\Om;\R^3),
\]
\[
 \mathbf{B}_0 \in C^2(\Om;\R^3),\,\, \Div_x \mathbf{B}_0=0.
\]

Then the initial value problem \eqref{lp2}--\eqref{lp4} admits infinitely many weak solutions in $(0,T)\times \Om$ emanating from the same initial data.
\end{Theorem}
\begin{Remark}
In fact, the weak solutions obtained in Theorem \ref{eom1} are more regular except for the velocity field and magnetic field. More precisely, the equation of continuity is satisfied in the strong sense, while \eqref{lp2}$_2$ and \eqref{lp2}$_3$ are satisfied in the weak sense, cf. relations \eqref{lm3}, \eqref{lm6}, and \eqref{lm5} below.
\end{Remark}

The following two sections are devoted to the proof of Theorem \ref{eom1}.
The arguments are adaptations of the method developed by De Lellis and Sz\'ekelyhidi \cite{DS1,DS2} to compressible setting in the spirit of \cite{CFK}.
We refer to \cite{C1,DFM,F1,FGS} for similar results on other physical models of related to compressible fluids.
The MHD case, however, is more delicate than the situation studied in \cite{CFK}
as there is no obvious way how to control the $L^\infty$-norm of the magnetic field $\mathbf{B}$.

\section{Reformulation and subsolutions }\label{ref}
First we reformulate \eqref{lp2} in terms of the new variables $[\vr,\vc{m}:=\vr \vv,\mathbf{B}]$:
\begin{equation}\label{lm1}
\left\{\begin{aligned}
& \p_t\vr + \Div_x \vc{m}=0,\\
& \p_t \vc{m}+\Div_x\left(\f{\vc{m} \otimes \vc{m}}{\vr}-\mathbf{B}\otimes \mathbf{B}\right)+\Grad_x\left(p(\vr)+\f{1}{2}|\mathbf{B}|^2\right)=
\mathbf{0},\\
& \p_t \mathbf{B}=\mathbf{curl}_x \left(\f{\vc{m} \times \mathbf{B}}{\vr} \right)+\Delta \mathbf{B}.\\
\end{aligned}\right.
\end{equation}
To facilitate the use of convex integration, we follow the strategy of \cite{CFK} writing
\[
\vc{m}=\mathbf{w}+\Grad_x \Phi,\,\,\Div_x \mathbf{w}=0,\,\,\int_{\Om}\Phi \ \dx=0.
\]
Denoting the traceless part of $\mathbf{B}\otimes \mathbf{B}$ as
\[
\mathbf{B} \odot \mathbf{B}:=\mathbf{B}\otimes \mathbf{B}-\f{1}{3}|\mathbf{B}|^2 \mathbf{I}_3,
\]
we rewrite  (\ref{lm1}) in the form
\begin{equation}\label{lm2}
\left\{\begin{aligned}
& \p_t\vr + \Delta \Phi=0,\\
& \p_t \mathbf{w}+\Div_x\left(\f{(\mathbf{w}+\Grad_x \Phi) \otimes (\mathbf{w}+\Grad_x \Phi)} {\vr}-\mathbf{B} \odot \mathbf{B}\right)\\
& +\Grad_x\left(p(\vr)+\f{1}{6}|\mathbf{B}|^2+\p_t \Phi\right)=
\mathbf{0},\\
& \p_t \mathbf{B}=\mathbf{curl}_x \left(\f{(\mathbf{w}+\Grad_x \Phi)\times \mathbf{B}}{\vr} \right)+\Delta \mathbf{B}.\\
\end{aligned}\right.
\end{equation}

Next, we use the ansatz for the density already exploited in \cite{CFK,F1}; specifically, we fix the density $\vr$,
\begin{equation}\label{lm9}
\left\{\begin{aligned}
& \vr \in C^2([0,T]\times \Om),\,\, \vr(0,\cdot)=\vr_0,\\
& \vr(t,x)>0,  \text{ for any } (t,x)\in [0,T]\times \Om.\\
\end{aligned}\right.
\end{equation}
With the density $\vr$ fixed, the potential function $\Phi$ is then uniquely solved by the Poisson equation
\begin{equation}\label{lm3}
\left\{\begin{aligned}
& -\Delta \Phi (t,\cdot)=\p_t \vr(t,\cdot),\\
& \int_{\Om}\Phi(t,\cdot)\dx=0.\\
\end{aligned}\right.
\end{equation}
Thus the continuity equation (\ref{lp2})$_1$ is satisfied in the strong sense.

Next, we observe that the equation of magnetic field (\ref{lm2})$_3$ is linear with respect to $\mathbf{B}$ for any given $\mathbf{w}=(w_1,w_2,w_3)$. Furthermore, we write
\[
\p_t \mathbf{B}- \p_{x_j}\left( \mathcal{A}_{i,j}\p_{x_i}\mathbf{B}+\mathcal{C}_j \mathbf{B}  \right)=\mathbf{0},
\]
where we have adopted the summation convention over repeated indices ($i,j$=1,2,3), and where
\begin{equation*}
\mathcal{A}_{i,j}=
\left\{\begin{aligned}
& \mathbf{I}_3\,\, \text{ if }i=j,\\
& \mathbf{0}\,\,\text{ if }i\neq j,\\
\end{aligned}\right.
\end{equation*}
\[
\mathcal{C}_1=
\begin{pmatrix}
0&0&0\\
\f{w_2+\p_{x_2}\Phi}{\vr}&-\f{w_1+\p_{x_1}\Phi}{\vr}&0\\
\f{w_3+\p_{x_3}\Phi}{\vr}&0&-\f{w_1+\p_{x_1}\Phi}{\vr}
\end{pmatrix},
\]
\[
\mathcal{C}_2=
\begin{pmatrix}
-\f{w_2+\p_{x_2}\Phi}{\vr}&\f{w_1+\p_{x_1}\Phi}{\vr}&0\\
0&0&0\\
0&\f{w_3+\p_{x_3}\Phi}{\vr}&-\f{w_2+\p_{x_2}\Phi}{\vr}
\end{pmatrix},
\]
\[
\mathcal{C}_3=
\begin{pmatrix}
-\f{w_3+\p_{x_3}\Phi}{\vr}&0&\f{w_1+\p_{x_1}\Phi}{\vr}\\
0&-\f{w_3+\p_{x_3}\Phi}{\vr}&\f{w_2+\p_{x_2}\Phi}{\vr}\\
0&0&0
\end{pmatrix}.
\]

Therefore, for any given $\mathbf{w}\in L^{\infty}((0,T)\times \Om;\R^3)$, the well-known maximal regularity for parabolic equations (see \cite{Amann,Kry}) implies that there exists a weak solution $\mathbf{B}=\mathbf{B}[\mathbf{w}]$ to the initial value problem
\begin{equation}\label{lm4}
\left\{\begin{aligned}
& \p_t \mathbf{B}=\mathbf{curl}_x \left(\f{(\mathbf{w}+\Grad_x \Phi)\times \mathbf{B}}{\vr} \right)+\Delta \mathbf{B},\\
& \mathbf{B}(0,\cdot)=\mathbf{B}_0,\\
\end{aligned}\right.
\end{equation}
unique in the class
\begin{equation}\label{lm5}
\left\{\begin{aligned}
& \mathbf{B }\in L^p(0,T;W^{1,p}(\Om;\R^3)),\\
& \p_t \mathbf{B} \in L^p(0,T;W^{-1,p}(\Om;\R^3)),\\
& \mathbf{B}\in C([0,T];\{W^{-1,p};W^{1,p} \}_{\alpha}),
\end{aligned}\right.
\end{equation}
for any $p\in (1,\infty)$, with suitable $\alpha \in (0,1)$, where $\{W^{-1,p};W^{1,p} \}_{\alpha}$ is the real interpolation space between $W^{-1,p}(\Om;\R^3)$ and $W^{1,p}(\Om;\R^3)$. We remark that (\ref{lm4}) and the assumption $\Div_x \mathbf{B}_0=0$ yield $\Div_x \mathbf{B}=0$ for a.e. $t\in(0,T)$ in the sense of distributions.

Thus the mapping $\mathbf{w} \mapsto \mathbf{B}[\mathbf{w}]$ can be considered as an abstract operator, whereas
the bounds in (\ref{lm5}) depend on the $L^{\infty}((0,T)\times \Om;\R^3)$-norm of $\mathbf{w}$. Therefore, it remains to prove that the ``abstract'' momentum equation:
\begin{equation}\label{lm6}
\left\{\begin{aligned}
& \p_t \mathbf{w}+\Div_x\left(\f{(\mathbf{w}+\Grad_x \Phi) \otimes (\mathbf{w}+\Grad_x \Phi)} {\vr}-\mathbf{B}[\mathbf{w}] \odot \mathbf{B}[\mathbf{w}]\right)\\
& +\Grad_x\left(p(\vr)+\f{1}{6}|\mathbf{B}[\mathbf{w}]|^2+\p_t \Phi-\f{2}{3}\eta\right)=
\mathbf{0},\\
& \Div_x \mathbf{w}=0,\,\,\mathbf{w}(0,\cdot)=\mathbf{w}_0:=\vr_0\vv_0-\Grad_x \Phi(0,\cdot).\\
\end{aligned}\right.
\end{equation}
admits infinitely many weak solutions in $(0,T)\times \Om$.
Similarly to \cite{CFK,F2}, we added a spatially homogeneous function $\eta=\eta(t),t\in [0,T]$, which is useful when
adjusting suitable energy bounds in the method of convex integration specified below.

In the text below we use the following notation. The symbol $\R^{N\times N}_{sym}$ stands for the space of $N\times N$ symmetric matrices over $\R$; $\R^{N\times N}_{sym,0}$ is the subspace of $\R^{N\times N}_{sym}$ with vanishing trace. Given $\mathbf{U}\in \R^{N\times N}_{sym}$, we denote by $\lambda_{max}[\mathbf{U}]$ its maximum eigenvalue. $C_w([0,T];L^2(\Om;\R^3))$ denotes the space of continuous functions from $[0,T]$ to $L^2(\Om;\R^3)$ equipped with the weak topology.

We are now in a position to introduce the class of subsolutions.
Inspired by \cite{F2}, we set the kinetic energy $e[\mathbf{w}]$ as
\beq\label{lm7}
e[\mathbf{w}]:=\eta(t)-\f{3}{2}p(\vr)-\f{1}{4}|\mathbf{B}[\mathbf{w}]|^2-\f{3}{2}\p_t \Phi
\eeq
for any $\mathbf{w}\in L^{\infty}((0,T)\times \Om;\R^3)$.
\begin{Remark}
The structure of the MHD system plays a crucial role in the definition of the kinetic energy $e[\mathbf{w}]$. Specifically, the
minus sign in the term ``$-\f{1}{4}|\mathbf{B}[\mathbf{w}]|^2 $" implies boundedness of the set of subsolutions
as detailed below.
\end{Remark}

In analogy with \cite{CFK,F2}, we define the space of subsolutions as
\[
X_0:=\Big\{\mathbf{w}\,\,|\mathbf{w}\in L^{\infty}((0,T)\times \Om;\R^3) \cap C_w([0,T];L^2(\Om;\R^3))
\cap C^1((0,T)\times \Om;\R^3),
\]
\[
\mathbf{w} \text{ satisfies the linear system }
\]
\begin{equation*}
\left\{\begin{aligned}
& \p_t \mathbf{w}+\Div_x \mathbf{U}=\mathbf{0},\\
& \Div_x \mathbf{w}=0,\,\,\mathbf{w}(0,\cdot)=\mathbf{w}_0,\\
\end{aligned}\right.
\end{equation*}
\centerline{for some $\mathbf{U}\in C^1((0,T)\times \Om;\R^{3\times 3}_{sym,0})\cap L^{\infty}((0,T)\times \Om;\R^{3\times 3}_{sym,0})$,}
\[
\f{3}{2}\lambda_{max}\left[\f{(\mathbf{w}+\Grad_x \Phi) \otimes (\mathbf{w}+\Grad_x \Phi)} {\vr}-\mathbf{B}[\mathbf{w}] \odot \mathbf{B}[\mathbf{w}]-\mathbf{U}  \right] < e[\mathbf{w}] \text{ in }(0,T)\times \Om \Big\}.
\]
Observe that we can choose $\eta\in C([0,T])$ in such a way that
\[
\f{3}{2}\lambda_{max}\left[\f{(\mathbf{w}_0+\Grad_x \Phi) \otimes (\mathbf{w}_0+\Grad_x \Phi)} {\vr}-\mathbf{B}[\mathbf{w_0}] \odot \mathbf{B}[\mathbf{w}_0] \right]
\]
\beq\label{lm8}
< \eta(t)-\f{3}{2}p(\vr)-\f{1}{4}|\mathbf{B}[\mathbf{w}_0]|^2-\f{3}{2}\p_t \Phi
\eeq
for any $t\in [0,T],x\in \Om$. This means that $\mathbf{w}_0$, together with the associated matrix field $\mathbf{U}=\mathbf{0}$ belongs to the space of subsolutions. In other words, $X_0$ is non-empty.

\section{Oscillatory lemma and convex integration scheme}\label{osc}

Following \cite{CFK}, we adapt the convex integration method from \cite{DS1} to the present setting.

\subsection{Oscillatory lemma}
The building block for the technique of convex integration developed by De Lellis and Sz\'ekelyhidi is the oscillatory lemma \cite{DS1,DS2} in the context of incompressible Euler system. The oscillatory lemma was later adapted to the case of compressible flows by Chiodaroli \cite{C1}. Here, we report a generalized version from \cite{DFM}.
\begin{Lemma}\label{invl}
Assume that
\[
\mathbf{h}\in C((0,T)\times \Om;\R^3),\,\, \mathbf{H}\in C((0,T)\times \Om; \R^{3\times 3}_{sym,0}),
\]
\[
e\in C((0,T)\times \Om),\,\, r \in C((0,T)\times \Om),\,\,r>0,\,\,0\leq e< \overline{e} \text{ in }(0,T)\times \Om.
\]
Assume also that
\[
\f{3}{2}\lambda_{max}\left[   \f{\mathbf{h}\otimes \mathbf{h}}{r}-\mathbf{H}
\right]< e \text{ in }(0,T)\times \Om.
\]

Then there exist two sequences
\[
\{ \mathbf{h}_n \}_{n\geq 1}\subset C_c^{\infty}((0,T)\times \Om;\R^3),\,\, \{ \mathbf{H}_n \}_{n\geq 1} \subset
C_c^{\infty}((0,T)\times \Om;\R^{3\times 3}_{sym,0})
\]
such that
\begin{equation}\label{jn6}
\left\{\begin{aligned}
& \p_t \mathbf{h}_n +\Div_x \mathbf{H}_n=0,\,\, \Div_x \mathbf{h}_n=0 \text{ in }(0,T)\times \Om,\\
& \f{3}{2} \lambda_{max}\left[   \f{(\mathbf{h}+\mathbf{h}_n)\otimes (\mathbf{h}+\mathbf{h}_n)}{r}-(\mathbf{H}+\mathbf{H}_n)
\right]< e \text{ in }(0,T)\times \Om,\\
& \mathbf{h}_n \rightarrow \mathbf{0} \text{ in } C_w([0,T];L^2(\Om;\R^3)),\\
& \liminf_{n\rightarrow \infty}\int_0^T \int_{\Om}\f{|\mathbf{h}_n|^2}{r}\dxdt
\geq C(\overline{e})\int_0^T\int_{\Om} \left( e-\f{1}{2}\f{|\mathbf{h}|^2}{r} \right)^2 \dxdt,\\
\end{aligned}\right.
\end{equation}
for some positive constant $C(\overline{e})$ depending only on $\overline{e}$.
\end{Lemma}

\subsection{Convex integration scheme}\label{consc}

We start with introducing the functional setting in which infinitely many weak solutions to (\ref{lm6}) will be located. To this end, observe that, as proved by De Lellis and Sz\'ekelyhidi \cite{DS2},
\beq\label{jn7}
\f{1}{2}|\mathbf{w}|^2\leq \f{3}{2}\lambda_{max}(\mathbf{w} \otimes \mathbf{w}-\mathbf{U})
\eeq
for any $\mathbf{w}\in \R^3,\,\mathbf{U}\in \R^{3\times 3}_{sym,0}$. Moreover, the equality holds if and only if
\beq\label{jn8}
\mathbf{U}=\mathbf{w} \otimes \mathbf{w}-\f{1}{3}|\mathbf{w}|^2 \mathbf{I}_3.
\eeq
From (\ref{lm9}), (\ref{lm7}), (\ref{jn7}) and the definition of the space $X_0$, we deduce the uniform bound
\beq\label{jn9}
|\mathbf{w}(t,x)|\leq \overline{w},\,\, t\in (0,T),\,x\in \Om
\eeq
for any $\mathbf{w}\in X_0$. This shows that the set of subsolutions $X_0$ is bounded in $L^{\infty}((0,T)\times \Om;\R^3)$. As an immediate consequence, we conclude that the bounds in (\ref{lm5}) are uniform with respect to $\mathbf{w}\in X_0$.
 Notice that (\ref{jn9}) implies that the functions belonging to $X_0$ range in a bounded ball $\Xi$ of $L^{2}( \Om;\R^3)$ for any $t\in [0,T]$, making the weak topology metrizable. The metric in $\Xi$ naturally induces a metric in $C([0,T];\Xi)$ denoted by $d$. We then define $X$ to be the closure of $X_0$ in $C_w([0,T];L^2(\Om;\R^3))$ with respect to the metric $d$. As a consequence, $(X,d)$ becomes a complete metric space bounded in $L^{\infty}((0,T)\times \Om;\R^3)$.

Following \cite{DS2,CFK,F2}, we introduce the functional $\mathcal{I} $ from $X$ to $(-\infty,0]$ as
\[
\mathcal{I}[\mathbf{w}]:=\int_0^T\int_{\Om}\left(\f{1}{2}\f{|\mathbf{w}+\Grad_x \Phi|^2}{\vr}
-e[\mathbf{w}] \right)\dxdt,\,\,\,\mathbf{w} \in X.
\]
It follows that the image of $\mathcal{I}$ ranges in a bounded interval of $(-\infty,0]$.

Next, we show strong continuity of the operators $\mathbf{w} \mapsto \mathbf{B}[\mathbf{w}]$. To this end, suppose
\beq\label{jn91}
\mathbf{w}_n\rightarrow \mathbf{w} \text{ in } X.
\eeq
On the one hand, we notice that interpolation space $\{W^{-1,p};W^{1,p} \}_{\alpha}$
appearing in the maximal regularity estimates is compactly embedded into $C(\Om;\R^3)$ for $p$ large enough, cf. \cite{Amann}. On the other hand, (\ref{lm5})$_2$ and (\ref{lm5})$_3$ together give the compactness in time. Thus we obtain
\beq\label{jn92}
\mathbf{B}[\mathbf{w}_n]\rightarrow \mathbf{B}[\mathbf{w}] \text{ in } C([0,T]\times \Om);
\eeq
whence
\beq\label{jn93}
e[\mathbf{w}_n]\rightarrow e[\mathbf{w}] \text{ in } C([0,T]\times \Om).
\eeq
The functional $\mathcal{I}$ is thus lower-semicontinuous on $X$ and belongs to Baire-1 mapping. Consequently, the points of continuity of $\mathcal{I}$ form a residual set in $X$ due to the well-known Baire's category theorem. To proceed, similarly to \cite{CFK,DS2,F2}, we employ the following crucial claim, the proof of which leans on the oscillatory lemma \ref{invl} and will be postponed to the end of this section.

{\bf{Claim.}} If $\mathbf{w}$ is a point of continuity of $\mathcal{I}$ in $X$, then $\mathcal{I}[\mathbf{w}]=0$.

Due to the oscillatory lemma \ref{invl} and the fact that $X_0$ is non-empty, the set
\[
\mathcal{M}:=\left\{  \mathbf{w}\in X\,\,|\,\,\mathcal{I} \text{ is continuous at }\mathbf{w}
\right\}
\]
admits infinite cardinality. Upon applying the claim above, we find that for any $\mathbf{w}\in \mathcal{M}$
\beq\label{jn10}
\f{1}{2}\f{|\mathbf{w}+\Grad_x \Phi|^2}{\vr}
=e[\mathbf{w}]=\eta(t)-\f{3}{2}p(\vr)-\f{1}{4}|\mathbf{B}[\mathbf{w}]|^2-\f{3}{2}\p_t \Phi
\eeq
for a.e. $t\in (0,T),x\in \Om$. In addition, it is also clear that $\mathbf{w}(0,\cdot)=\mathbf{w}_0$ and
\begin{equation*}
\left\{\begin{aligned}
& \p_t \mathbf{w}+\Div_x \mathbf{U}=\mathbf{0},\\
& \Div_x \mathbf{w}=0,\\
\end{aligned}\right.
\end{equation*}
in the sense of distributions. Recalling the definition of subsolutions $X_0$ and (\ref{jn7})-(\ref{jn8}), we arrive at
\[
\mathbf{U}=\f{(\mathbf{w}+\Grad_x \Phi) \otimes (\mathbf{w}+\Grad_x \Phi)} {\vr}-\f{1}{3}\f{|\mathbf{w}+\Grad_x \Phi|^2}{\vr}\mathbf{I}_3
\]
for a.e. $t\in (0,T),x\in \Om$. This shows every $\mathbf{w}$ in $\mathcal{M}$ exactly solves (\ref{lm6}). Therefore, it remains to verify the claim above so as to finish the proof of Theorem \ref{eom1}.

For completeness, we provide the detailed proof of the claim following the arguments of \cite{F2} in the context of an abstract Euler-type system. Let $\mathbf{w}$ be a point of continuity of $\mathcal{I}$ in $X$. Then there exists a sequence $\{\mathbf{w}_k\}_{k \geq 1}\subset X_0$ such that as $k\rightarrow \infty$
\[
\mathbf{w}_k \rightarrow \mathbf{w} \text{ in }C_w([0,T];L^2(\Om;\R^3)).
\]
By the assumption of continuity of $\mathcal{I}$,
\[
\mathcal{I}[\mathbf{w}_k]\rightarrow \mathcal{I}[\mathbf{w}] \text{ as }k\rightarrow \infty.
\]
In accordance with the definition of subsolutions $X_0$, there exists a suitable sequence $\{\ep_{k}\}_{k \geq 1}$ tending to zero such that
\[
\f{3}{2}\lambda_{max}\left[\f{(\mathbf{w}_k+\Grad_x \Phi) \otimes (\mathbf{w}_k+\Grad_x \Phi)} {\vr}-\mathbf{B}[\mathbf{w}_k] \odot \mathbf{B}[\mathbf{w}_k]-\mathbf{U}_k  \right]
\]
\[
< e[\mathbf{w}_k]-\ep_k \text{ in }(0,T)\times \Om,
\]
where $\{\mathbf{U}_k\}_{k \geq 1} $ are the fluxes associated with $\{\mathbf{w}_k\}_{k \geq 1}$.
At this stage we fix $k$ and apply the oscillatory lemma \ref{invl} with
\[
r=\vr,\,\,\mathbf{h}=\mathbf{w}_k+\Grad_x \Phi,\,\,\mathbf{H}=\mathbf{B}[\mathbf{w}_k] \odot \mathbf{B}[\mathbf{w}_k]+\mathbf{U}_k,\,\,e=e[\mathbf{w}_k]-\ep_k,
\]
to deduce that there exist compactly supported smooth sequences $\{\mathbf{h}_k^n\}_{n\geq 1},\,\{ \mathbf{H}_k^n \}_{n\geq 1}$ obeying the properties therein. In particular, by setting the new variables
\[
\mathbf{w}_k^n:=\mathbf{w}_k+\mathbf{h}_k^n,\,\,\mathbf{U}_k^n:=\mathbf{U}_k+\mathbf{H}_k^n,
\]
we see
\begin{equation*}
\left\{\begin{aligned}
& \p_t \mathbf{w}_k^n+\Div_x \mathbf{U}_k^n=\mathbf{0},\\
& \Div_x \mathbf{w}_k^n=0,\,\,\mathbf{w}_k^n(0,\cdot)=\mathbf{w}_0,\\
\end{aligned}\right.
\end{equation*}
\[
\f{3}{2}\lambda_{max}\left[\f{(\mathbf{w}_k^n+\Grad_x \Phi) \otimes (\mathbf{w}_k^n+\Grad_x \Phi)} {\vr}-\mathbf{B}[\mathbf{w}_k] \odot \mathbf{B}[\mathbf{w}_k]-\mathbf{U}_k^n  \right]
\]
\[
< e[\mathbf{w}_k]-\ep_k \text{ in }(0,T)\times \Om.
\]
Applying the continuity of $e[\cdot]$ (see (\ref{jn91}), (\ref{jn93})), we conclude
\[
e[\mathbf{w}_k^n]\rightarrow e[\mathbf{w}_k] \text{ in }C([0,T]\times \Om)
\]
as $n\rightarrow \infty$. Consequently, for each $k\geq 1$, there exists $n=n(k)$ such that
\[
\mathbf{w}_k^{n(k)}\in X_0.
\]
Obviously, in light of (\ref{jn6})$_3$,
\[
\mathbf{w}_k^{n(k)} \rightarrow \mathbf{w} \text{ in }C_w([0,T];L^2(\Om;\R^3))
\]
as $k\rightarrow \infty$. The continuity of $\mathcal{I}$ again implies
\[
\mathcal{I}[\mathbf{w}_k^{n(k)} ]\rightarrow \mathcal{I}[\mathbf{w}] \text{ as }k\rightarrow \infty.
\]
Finally, by virtue of (\ref{jn6})$_3$, (\ref{jn6})$_4$ and Cauchy-Schwarz's inequality, we obtain
\[
\liminf_{k\rightarrow \infty}\mathcal{I}[\mathbf{w}_k^{n(k)}]
=\liminf_{k\rightarrow \infty}\int_0^T\int_{\Om}
\left(\f{1}{2}\f{|\mathbf{w}_{k}+\mathbf{h}_k^{n(k)}+\Grad_x \Phi|^2}{\vr}
-e\left[\mathbf{w}_k+\mathbf{h}_k^{n(k)}\right] \right)\dxdt
\]
\[
=\lim_{k\rightarrow \infty}\int_0^T\int_{\Om}
\left(\f{1}{2}\f{|\mathbf{w}_{k}+\Grad_x \Phi|^2}{\vr}
-e\left[\mathbf{w}_k+\mathbf{h}_k^{n(k)}\right] \right)\dxdt
+\liminf_{k\rightarrow \infty} \int_0^T\int_{\Om} \f{1}{2}\f{|\mathbf{h}_k^{n(k)}|^2}{\vr} \dxdt
\]
\[
\geq \mathcal{I}[\mathbf{w}]+\f{C(\overline{e})}{2}\liminf_{k\rightarrow \infty}\int_0^T\int_{\Om}
\left( e[\mathbf{w}_k]-\ep_k- \f{1}{2} \f{|\mathbf{w}_{k}+\Grad_x \Phi|^2}{\vr}  \right)^2
\dxdt
\]
\[
\geq \mathcal{I}[\mathbf{w}]+\f{C(\overline{e})}{2T}\liminf_{k\rightarrow \infty}\left[\int_0^T\int_{\Om}
\left( e[\mathbf{w}_k]-\ep_k- \f{1}{2} \f{|\mathbf{w}_{k}+\Grad_x \Phi|^2}{\vr}  \right)
\dxdt\right]^2
\]
\[
=\mathcal{I}[\mathbf{w}]+\f{C(\overline{e})}{2T}\left(  \mathcal{I}[\mathbf{w}] \right)^2.
\]
We conclude from the inequality above that $\mathcal{I}[\mathbf{w}]=0$, thus verifying the claim.

We have completed the proof of Theorem \ref{eom1}.

\section{Admissible solutions}\label{conrm}

Let us denote by $H(\vr)$ the potential energy of (\ref{lp2})
\beq\label{lu00}
H(\vr):=\vr\int_1^{\vr}\f{p(s)}{s^2}\ds;
\eeq
$E(t)$ the total energy
\beq\label{lu01}
E(t):=\int_{\Om} \left( \f{1}{2}\vr |\vv|^2+\f{1}{2}|\mathbf{B}|^2+H(\vr) \right)(t,x)\dx.
\eeq
Assume for the moment that $[\vr,\vv,\mathbf{B}]$ is a smooth solution to the problem (\ref{lp2})-(\ref{lp4}). It is easy to derive the energy inequality
\[
\f{d}{dt}E(t)+\int_{\Om}|\Grad_x \mathbf{B}|^2\dx\leq0.
\]
This particularly implies
\beq\label{lu1}
E(t)\leq E(0) \text{ for a.e. } t\in (0,T).
\eeq
However, the relation (\ref{lu1}) is no longer available for the weak solutions obtained in Theorem \ref{eom1} by the technique of convex integration. As a matter of fact, it can be seen from (\ref{lm8}) and (\ref{jn10}) that
\[
\liminf_{t\rightarrow 0^+}E(t)>E(0).
\]
This is the initial energy jump typical for weak solutions obtained through the method of convex integration; see \cite{CM1,CFK,F1}. Nevertheless, it is possible to remove this drawback by modifying the convex integration scheme as explained in \cite{CFK,F2}, at least for certain initial velocity.
\begin{Theorem}\label{eom2}
Let $T\in (0,\infty)$ be given. Suppose that
\[
p(\cdot)\in C([0,\infty))\cap C^\infty((0,\infty)),\,\, p(0)=0,
\]
\[
\vr_0\in C^3(\Om),\,\,\vr_0>\underline{\vr}>0,\,\, \mathbf{B}_0 \in C^2(\Om;\R^3),\,\,\Div_x \mathbf{B}_0=0.
\]

Then there exists $\vv_0\in L^{\infty}(\Om;\R^3)$ such that the initial value problem (\ref{lp2})-(\ref{lp4}) admits infinitely many weak solutions in $(0,T)\times \Om$ satisfying the energy inequality (\ref{lu1}).
\end{Theorem}

The proof of Theorem \ref{eom2} follows the same lines as \cite{CFK}, using the convex integration scheme given in Sections \ref{ref}-\ref{osc}. We thus omit the details here.

\section{Extension to heat-conductive fluid flows}\label{exth}
In this section, we show the existence and non-uniqueness of global-in-time weak solutions to the compressible MHD system for inviscid resistive and heat conductive fluids. By assuming the fluid to be ideal and polytropic, the governing equations take the form
\begin{equation}\label{jn11}
\left\{\begin{aligned}
& \p_t\vr + \Div_x(\vr\vv)=0,\\
& \p_t(\vr\vv)+\Div_x\left(\vr\vv \otimes \vv-\mathbf{B}\otimes \mathbf{B}\right)+\Grad_x\left( \vr\vt+\f{1}{2}|\mathbf{B}|^2\right)=
\mathbf{0},\\
& \p_t \mathbf{B}=\mathbf{curl}_x (\vv \times \mathbf{B})+\Delta \mathbf{B},\\
& \Div_x \mathbf{B}=0,\\
& c_V \left[ \p_t (\vr\vt)+\Div_x (\vr\vt\vv) \right]=\Delta \vt - \vr \vt \Div_x \vv +
|\mathbf{curl}_x \mathbf{B}|^2. \\
\end{aligned}\right.
\end{equation}
The unknowns $\vr,\vv$ and $\mathbf{B}$ are defined as before, $\vt$ is the (absolute) temperature, and $c_V>0$ denotes the specific heat at constant volume. For technical reason (see Remark \ref{mrk1}), we are not able to handle the general MHD system (\ref{jn11}). Inspired by \cite{BLN,LIS}, we shall consider the three-dimensional system with certain symmetry. More precisely, by setting
\begin{equation}\label{lp5}
\left\{\begin{aligned}
& \vv(x,t)=(u_1(x_1,x_2,t),u_2(x_1,x_2,t),0)=:(\vu,0),\\
& \vr(x,t)=\vr(x_1,x_2,t),\\
& \vt(x,t)=\vt(x_1,x_2,t),\\
& \mathbf{B}(x,t)=(0,0,b(x_1,x_2,t)),\\
\end{aligned}\right.
\end{equation}
in (\ref{jn11}), one obtains the following 2-D system of equations:
\begin{equation}\label{lp3}
\left\{\begin{aligned}
& \p_t\vr + \Div_x(\vr\vu)=0,\\
& \p_t(\vr\vu)+\Div_x(\vr\vu \otimes \vu)+\Grad_x\left( \vr\vt+\f{1}{2}b^2\right)=
\mathbf{0},\\
& \p_t b+ \Div_x(b\vu)=\Delta b,\\
& c_V \left[ \p_t (\vr\vt)+\Div_x (\vr\vt\vu)          \right]=\Delta \vt - \vr \vt \Div_x \vu +
|\Grad_x b|^2. \\
\end{aligned}\right.
\end{equation}
Here we use the same symbol $x$ to denote the two-dimensional space variable $(x_1,x_2)$ for convenience. Similarly to the above, we consider the periodic boundary conditions
\[
\Om:= \left( [0,1]|_{\{0,1\}} \right)^2,
\]
and (\ref{lp3}) is supplemented with the initial conditions:
\beq\label{jn121}
\vr(0,\cdot)=\vr_0,\,\, \vu(0,\cdot)=\vu_0,\,\,b(0,\cdot)=b_0,\,\,\vt(0,\cdot)=\vt_0.
\eeq
The definition of weak solutions to the problem (\ref{lp3})-(\ref{jn121}) is analogous to the general barotropic case (\ref{lp2}). Then our result of non-uniqueness of global-in-time weak solutions to (\ref{lp3})-(\ref{jn121}) can be stated as follows.
\begin{Theorem}\label{eom4}
Let $T\in (0,\infty)$ be given. Assume that
\[
p(\cdot)\in C([0,\infty))\cap C^\infty((0,\infty)),\,\, p(0)=0,
\]
\[
\vr_0\in C^3(\Om),\,\,\vr_0>\underline{\vr}>0,\,\,\vu_0 \in C^3(\Om;\R^2),
\]
\[
b_0 \in C^2(\Om),\,\,\vt_0\in C^2(\Om), \vt_0>\underline{\vt}>0.
\]

Then the initial value problem (\ref{lp3})-(\ref{jn121}) admits infinitely many weak solutions in $(0,T)\times \Om$ starting from the same initial data.
\end{Theorem}
\begin{Remark}
In contrast with the general barotropic case, the weak solutions obtained in Theorem \ref{eom4} satisfy the equation of continuity,
the equations for the {magnetic field} and the temperature in the strong sense, while the momentum equation is satisfied in the weak sense.
\end{Remark}

\begin{Remark}
Since every weak solution to the symmetry reduced system (\ref{lp3}) can be regarded as a weak solution to the three-dimensional system (\ref{jn11}) through the identification (\ref{lp5}), Theorem \ref{eom4} in particular provides infinitely many global-in-time weak solutions to (\ref{jn11}).
\end{Remark}

The proof of Theorem \ref{eom4} follows basically the same lines as Theorem \ref{eom1} with only slight modifications. We thus give the outline of proof. In the new variables $[\vr,\vc{m}=\vr \vu,b,\vt]$, (\ref{lp3}) is reformulated as
\begin{equation}\label{jn12}
\left\{\begin{aligned}
& \p_t\vr + \Div_x \vc{m}=0,\\
& \p_t \vc{m}+\Div_x\left(\f{\vc{m} \otimes \vc{m}}{\vr}\right)+\Grad_x\left(\vr \vt+\f{1}{2}b^2\right)=
\mathbf{0},\\
& \p_t b+ \f{b}{\vr}\Div_x \vc{m}-\f{b}{\vr^2}\vc{m}\cdot \Grad_x \vr+\f{1}{\vr}\vc{m}\cdot
 \Grad_x b= \Delta b,\\
& c_V [\vr \p_t \vt+\vc{m}\cdot \Grad_x \vt]=\Delta \vt -\vt \Div_x \vc{m}+\vt\f{\Grad_x \vr}{\vr}\cdot \vc{m}+|\Grad_x b|^2.\\
\end{aligned}\right.
\end{equation}
As before, we invoke the Helmholtz decomposition
\[
\vc{m}=\mathbf{w}+\Grad_x \Phi,\,\,\Div_x \mathbf{w}=0,\,\,\int_{\Om}\Phi \dx=0.
\]
Consequently, (\ref{jn12}) is rewritten as
\begin{equation}\label{jn13}
\left\{\begin{aligned}
& \p_t\vr + \Delta \Phi=0,\\
& \p_t \mathbf{w}+\Div_x\left(\f{(\mathbf{w}+\Grad_x \Phi) \otimes (\mathbf{w}+\Grad_x \Phi)} {\vr}\right)+\Grad_x\left(\vr\vt+\f{1}{2}b^2+\p_t \Phi\right)=
\mathbf{0},\\
& \p_t b+ \f{b}{\vr} \Delta \Phi-\f{b}{\vr^2}(\mathbf{w}+\Grad_x \Phi)\cdot \Grad_x \vr+\f{1}{\vr}(\mathbf{w}+\Grad_x \Phi)\cdot
 \Grad_x b= \Delta b,\\
& c_V [\vr \p_t \vt +(\mathbf{w}+\Grad_x \Phi)\cdot \Grad_x \vt]=\Delta \vt -\vt \Delta \Phi +\vt \f{\Grad_x \vr}{\vr}\cdot (\mathbf{w}+\Grad_x \Phi)+|\Grad_x b|^2.\\
\end{aligned}\right.
\end{equation}

Next, we adopt the ansatz for the density as in (\ref{lm9}). To proceed, we observe that for any given $\mathbf{w}\in L^{\infty}((0,T)\times \Om;\R^2)$, (\ref{jn13})$_3$ is linear with respect to $b$ and the standard $L^p$-theory for parabolic equations (see \cite{Amann,Kry}) shows that there exists a solution $b=b[\mathbf{w}]$ to
\begin{equation}\label{jn14}
\left\{\begin{aligned}
& \p_t b+ \f{b}{\vr} \Delta \Phi-\f{b}{\vr^2}(\mathbf{w}+\Grad_x \Phi)\cdot \Grad_x \vr+\f{1}{\vr}(\mathbf{w}+\Grad_x \Phi)\cdot
 \Grad_x b= \Delta b,\\
& b(0,\cdot)=b_0,\\
\end{aligned}\right.
\end{equation}
unique in the class
\beq\label{jn15}
\p_t b \in L^p(0,T;L^p(\Om)),\,\,\Grad_x^2 b \in L^p(0,T;L^p(\Om;\R^{2\times 2}))
\eeq
for any $p \in (1,\infty)$. Analogously, there exists a solution $\vt=\vt[\mathbf{w}]$ to
\begin{equation}\label{jn16}
\left\{\begin{aligned}
&  c_V [\vr \p_t \vt +(\mathbf{w}+\Grad_x \Phi)\cdot \Grad_x \vt]=\Delta \vt -\vt \Delta \Phi +\vt \f{\Grad_x \vr}{\vr}\cdot (\mathbf{w}+\Grad_x \Phi)+|\Grad_x b[\mathbf{w}]|^2,\\
& \vt(0,\cdot)=\vt_0,\\
\end{aligned}\right.
\end{equation}
unique in the class
\beq\label{jn17}
\p_t \vt \in L^p(0,T;L^p(\Om)),\,\,\Grad_x^2 \vt \in L^p(0,T;L^p(\Om;\R^{2\times 2}))
\eeq
for any $p \in (1,\infty)$. Therefore, the proof of Theorem \ref{eom4} reduces to find infinitely many solutions to
\begin{equation}\label{jn18}
\left\{\begin{aligned}
& \p_t \mathbf{w}+\Div_x\left(\f{(\mathbf{w}+\Grad_x \Phi) \otimes (\mathbf{w}+\Grad_x \Phi)} {\vr}\right)\\
& +\Grad_x\left(\vr\vt[\mathbf{w}]+\f{1}{2}b[\mathbf{w}]^2+\p_t \Phi-\eta(t)\right)=
\mathbf{0},\\
& \Div_x \mathbf{w}=0,\,\,\mathbf{w}(0,\cdot)=\mathbf{w}_0:=\vr_0\vu_0-\Grad_x \Phi(0,\cdot),\\
\end{aligned}\right.
\end{equation}
belonging to
\beq\label{jn19}
 L^{\infty}((0,T)\times \Om;\R^2) \cap C_w([0,T];L^2(\Om;\R^2)).
\eeq

To introduce the definition of subsolutions, we first modify the kinetic energy as
\[
e[\mathbf{w}]:=\eta(t)-\vr\vt[\mathbf{w}]-\f{1}{2}b[\mathbf{w}]^2-\p_t \Phi
\]
for any $\mathbf{w}\in L^{\infty}((0,T)\times \Om;\R^2)$. Then we define
\[
X_0:=\Big\{\mathbf{w}\,\,|\mathbf{w}\in L^{\infty}((0,T)\times \Om;\R^2) \cap C_w([0,T];L^2(\Om;\R^2))
\cap C^1((0,T)\times \Om;\R^2),
\]
\[
\mathbf{w} \text{ satisfies the linear system }
\]
\begin{equation*}
\left\{\begin{aligned}
& \p_t \mathbf{w}+\Div_x \mathbf{U}=\mathbf{0},\\
& \Div_x \mathbf{w}=0,\,\,\mathbf{w}(0,\cdot)=\mathbf{w}_0,\\
\end{aligned}\right.
\end{equation*}
\centerline{for some $\mathbf{U}\in C^1((0,T)\times \Om;\R^{2\times 2}_{sym,0})\cap L^{\infty}((0,T)\times \Om;\R^{2\times 2}_{sym,0})$,}
\[
\lambda_{max}\left[\f{(\mathbf{w}+\Grad_x \Phi) \otimes (\mathbf{w}+\Grad_x \Phi)} {\vr}-\mathbf{U}  \right] < e[\mathbf{w}] \text{ in }(0,T)\times \Om \Big\}.
\]
The function $\eta\in C([0,T])$ is then determined such that
\beq\label{jn20}
\lambda_{max}\left[\f{(\mathbf{w}_0+\Grad_x \Phi) \otimes (\mathbf{w}_0+\Grad_x \Phi)} {\vr} \right] < \eta(t)-\vr\vt[\mathbf{w}_0]-\f{1}{2}b[\mathbf{w}_0]^2-\p_t \Phi
\eeq
for any $t\in [0,T],x\in \Om$. Therefore, $X_0$ is non-empty since obviously it contains $\mathbf{w}_0$. Moreover, from the definition of subsolutions, one infers the uniform bound
\beq\label{jn21}
|\mathbf{w}(t,x)|\leq \overline{w}, \text{ for any } t\in [0,T],x\in \Om,\mathbf{w}\in X_0.
\eeq
It follows from (\ref{jn21}) and the parabolic estimates (\ref{jn15}), (\ref{jn17}) that
\beq\label{jn22}
e[\mathbf{w}_n]\rightarrow e[\mathbf{w}] \text{ in }C([0,T]\times \Om)
\eeq
whenever
\[
\mathbf{w}_n\rightarrow \mathbf{w} \text{ in } C_w([0,T];L^2(\Om;\R^2))
\text{ and weakly}-\ast \text{ in }L^{\infty}((0,T)\times \Om;\R^2).
\]
Then we repeat step by step the convex integration scheme as for the general barotropic case. The details are omitted.

\begin{Remark}\label{mrk1}
For the general three-dimensional heat-conductive MHD system (\ref{jn11}), property (\ref{jn22}) is unavailable. Precisely, we are lacking in the compactness of $\mathbf{curl}_x \mathbf{B}$.
\end{Remark}

We finish this section by proving the uniform bound of $b[\mathbf{w}]$ with respect to $\mathbf{w}\in L^{\infty}((0,T)\times \Om;\R^2)$. It is of independent interest. Following \cite{CFK}, we basically employ the comparison principle of parabolic equations. Furthermore, the structure of the MHD system (\ref{lp3}) plays a crucial role. The main observations are as follows.
\[
\p_t \left(\f{b}{\vr}\right) +\vu\cdot \Grad_x \left(\f{b}{\vr}\right)
\]
\[
=\f{1}{\vr}\p_t b -\f{b}{\vr^2}\p_t \vr +\f{1}{\vr}\vu \cdot \Grad _x b -\f{b}{\vr^2}\vu \cdot \Grad _x \vr
\]
\[
=\f{1}{\vr}\left(\Delta b-b \Div_x \vu \right)-\f{b}{\vr^2}(-\vr \Div_x \vu)
\]
\beq\label{jn1}
=\f{1}{\vr}\Delta b,
\eeq
where we used (\ref{lp3})$_1$ and (\ref{lp3})$_3$ in the second equality of (\ref{jn1}). A straightforward calculation gives that
\[
\Delta \left(\f{b}{\vr}\right) =
\f{1}{\vr}\Delta b +b \Delta \left( \f{1}{\vr} \right)+ 2 \Grad_x b \cdot \Grad_x \left( \f{1}{\vr} \right).
\]
From this relation we find
\[
\p_t \left(\f{b}{\vr}\right) +\vu\cdot \Grad_x \left(\f{b}{\vr}\right)
\]
\beq\label{jn2}
=\Delta \left(\f{b}{\vr}\right)-b \Delta \left( \f{1}{\vr} \right)-2 \Grad_x b \cdot \Grad_x \left( \f{1}{\vr} \right).
\eeq
Noticing that
\[
\Grad_x b=\vr \f{1}{\vr} \Grad_x b
\]
\beq\label{jn3}
=\vr \left[  \Grad_x \left(\f{b}{\vr}\right)-b \Grad_x \left( \f{1}{\vr} \right) \right].
\eeq
Therefore, combining (\ref{jn2}) with (\ref{jn3}), after returning to the variables $[\vr,\mathbf{w},b]$, yields
\[
\p_t \left(\f{b}{\vr}\right) + \left[   \f{1}{\vr}(\mathbf{w}+\Grad_x \Phi)  +2 \vr \Grad_x \left(\f{1}{\vr}\right)               \right]\cdot \Grad_x \left(\f{b}{\vr}\right)
\]
\beq\label{jn4}
=\Delta \left(\f{b}{\vr}\right)- \left[ \vr \Delta \left( \f{1}{\vr} \right)-2\vr^2 \left|\Grad_x \left( \f{1}{\vr} \right)\right|^2
 \right] \f{b}{\vr}.
\eeq
This shows that (\ref{jn4}) is linear with respect to $b/\vr$. As a consequence, in light of (\ref{lm9}), we apply the comparison principle of parabolic equations to (\ref{jn4}) to conclude that there exists a positive constant $\overline{b}$ depending only on the initial data and $T$ such that
\beq\label{jn5}
|b[\mathbf{w}](t,x)|\leq \overline{b}, \text{ for any }t \in [0,T], x\in \Om.
\eeq
Here, the point is that the bound $\overline{b}$ is independent of $\mathbf{w}$ belonging to $L^{\infty}((0,T)\times \Om;\R^2)$.

\medskip

\centerline{\bf Acknowledgement}

\medskip

E.Feireisl acknowledges support of the project 18-12719S financed
by the Czech Science Foundation.

The research of Y. Li is partially supported by  Postgraduate Research and Practice Innovation Program of Jiangsu Province under grant number KYCX 18-0028 and China Scholarship Council; he is also indebted to the Institute of Mathematics of the Czech Academy of Sciences for the invitation and hospitality.



\begin{thebibliography}{100}

\bibitem{Amann}Amann, H.: Linear and Quasilinear Parabolic Problems. Vol. I. Abstract linear theory. Monographs in Mathematics, {\bf89}. Birkh\"{a}user Boston, Inc., Boston(1995)

\bibitem{BSS}Bardos, C., Sulem, C., Sulem, P.L.: Longtime dynamics of a conductive fluid in the presence of a strong magnetic field. {Trans. Amer. Math. Soc.} {\bf305}, 175-191(1988)

\bibitem{BLN}Bronzi, A.C., Lopes Filho, M.C., Nussenzveig Lopes, H.J.: Wild solutions for 2D incompressible ideal flow with passive tracer. {Commun. Math. Sci.} {\bf13}, 1333-1343(2015)

\bibitem{HC}Cabannes, H.: {Theoretical Magnetofluiddynamics}. Academic Press, New York(1970)

\bibitem{CKS}Caflisch, R.E., Klapper, I., Steele, G.: Remarks on singularities, dimension and energy dissipation for ideal hydrodynamics and MHD. {Comm. Math. Phys.} {\bf184}, 443-455(1997)


\bibitem{CW}Cao, C., Wu, J.: Global regularity for the 2D MHD equations with mixed partial dissipation and magnetic diffusion. {Adv. Math.} {\bf226}, 1803-1822(2011)

\bibitem{C1}Chiodaroli, E.: A counterexample to well-posedness of entropy solutions to the compressible Euler system. {J. Hyperbolic Differ. Equ.} {\bf 11}, 493-519(2014)

\bibitem{CM1}Chiodaroli, E., Mich\'alek, M.: Existence and non-uniqueness of global weak solutions to inviscid primitive and Boussinesq equations. {Comm. Math. Phys.} {\bf 353} 1201-1216(2017)

\bibitem{CFK}Chiodaroli, E., Feireisl, E., Kreml, O.: On the weak solutions to the equations of a compressible heat-conducting gas. {Ann. I. H. Poincar\'e.} {\bf 32}, 225-243(2015)

\bibitem{CFG}C\'ordoba, D.: Faraco, D., Gancedo, F.: Lack of uniqueness for weak solutions of the incompressible porous media equation. {Arch. Ration. Mech. Anal.} {\bf 200}, 725-746(2011)

\bibitem{Dai} Dai, M.: Non-uniqueness of Leray-Hopf weak solutions of the 3D Hall-MHD system.
{arXiv: 1812.11311v2}(2019)

\bibitem{DS1}De Lellis, C., Sz\'ekelyhidi  Jr., L.: The Euler equations as a differential inclusion. {Ann. of Math. (2)} {\bf 170}, 1417-1436(2009)

\bibitem{DS2}De Lellis, C., Sz\'ekelyhidi  Jr., L.: On admissibility criteria for weak solutions of the Euler equations. {Arch. Ration. Mech. Anal.} {\bf 195}, 225-260(2010)

\bibitem{DFM}Donatelli, D., Feireisl, E., Marcati, P.: Well/ill posedness for the Euler-Korteweg-Poisson system and related problems. {Comm. Partial Differential Equations.} {\bf40}, 1314-1335(2015)

\bibitem{DF}Ducomet, B., Feireisl, E.: The equations of magnetohydrodynamics: on the interaction between matter and radiation in the evolution of gaseous stars. {Comm. Math. Phys.} {\bf266}, 595-629(2006)

\bibitem{DL}Duvaut, G., Lions, J.L.: In\'equations en thermo\'elasticit\'e et magn\'etohydrodynamique. {Arch. Rational Mech. Anal.} {\bf46}, 241-279(1972)

\bibitem{FJN}Fan, J., Jiang, S., Nakamura, G.: Vanishing shear viscosity limit in the magnetohydrodynamic equations. {Comm. Math. Phys.} {\bf 270}, 691-708(2007)

\bibitem{FL}Faraco, D., Lindberg, S.: Magnetic helicity and subsolutions in ideal MHD. {arXiv: 1801.04896}(2018)

\bibitem{F1}Feireisl, E.: On weak solutions to a diffuse interface model of a binary mixture of compressible fluids. {Discrete Contin. Dyn. Syst. Ser. S.} {\bf 9}, 173-183(2016)

\bibitem{F2}Feireisl, E.: Weak solutions to problems involving inviscid fluids. Mathematical fluid dynamics, present and future. {Springer Proc. Math. Stat.} {\bf 183}, 377-399(2016)

\bibitem{FGS}Feireisl, E., Gwiazda, P., Swierczewska-Gwiazda, A.: On weak solutions to the 2D Savage-Hutter model of the motion of a gravity-driven avalanche flow. {Comm. Partial Differential Equations.} {\bf 41}, 759-773(2016)

\bibitem{FKKM}Feireisl, E., Klingenberg, C., Kreml, O., Markfelder, S.: On oscillatory solutions to the complete Euler system. {arXiv: 1710.10918}(2017)

\bibitem{GMS}Gwiazda, P., Mich\'alek, M., \'Swierczewska-Gwiazda, A.: A note on weak solutions of conservation laws and energy/entropy conservation. {Arch. Ration. Mech. Anal.} {\bf229}, 1223-1238(2018)

\bibitem{HW}Hu, X., Wang, D.: Global existence and large-time behavior of solutions to the three-dimensional equations of compressible magnetohydrodynamic flows. {Arch. Rational Mech. Anal.} {\bf 197}, 203-238(2010)

\bibitem{KM}Klingenberg, C., Markfelder, S.: Non-uniqueness of entropy-conserving solutions to the ideal compressible MHD equations. {arXiv: 1902.01446}(2019)

\bibitem{Kry}Krylov, N.V.: Parabolic equations with VMO coefficients in Sobolev spaces with mixed norms. {J. Funct. Anal.} {\bf250}, 521-558(2007)

\bibitem{LXZ2}Li, H., Xu, X., Zhang, J.: Global classical solutions to 3D compressible magnetohydrodynamic equations with large oscillations and vacuum. {SIAM J. Math. Anal.} {\bf45}, 1356-1387(2013)


\bibitem{LIS}Li, Y., Sun, Y.: Global weak solutions to a two-dimensional compressible MHD equations of viscous non-resistive fluids. {arXiv: 1811.05752}(2018)

\bibitem{LXZ} Lin, F., Xu, L., Zhang, P.: Global small solutions of 2-D incompressible MHD system. { J. Differential Equations.} {\bf 259}, 5440-5485(2015)

\bibitem{LXX}Luo, T., Xie, C., Xin, Z.: Non-uniqueness of admissible weak solutions to compressible Euler systems with source terms. {Adv. Math.} {\bf 291}, 542-583(2016)

\bibitem{PZZ}Pan, R., Zhou, Y., Zhu, Y.: Global classical solutions of three dimensional viscous MHD system without magnetic diffusion on periodic boxes. {Arch. Ration. Mech. Anal.} {\bf227}, 637-662(2018)

\bibitem{PG}Pu, X., Guo, B.: Global existence and convergence rates of smooth solutions for the full compressible MHD equations. {Z. Angew. Math. Phys.} {\bf64}, 519-538(2013)

\bibitem{Sche}Scheffer, V.: An inviscid flow with compact support in space-time. {J. Geom. Anal.} {\bf 3}, 343-401(1993)

\bibitem{ST}Sermange, M., Temam, R.: Some mathematical questions related to the MHD equations. { Comm. Pure Appl. Math.} {\bf36}, 635-664(1983)

\bibitem{Shir}Shnirelman, A.: On the nonuniqueness of weak solution of the Euler equation. {Comm. Pure Appl. Math.} {\bf 50}, 1261-1286(1997)

\bibitem{Shvy}Shvydkoy, R.: Convex integration for a class of active scalar equations. {J. Amer. Math. Soc.} {\bf 24}, 1159-1174(2011)

\bibitem{TW} Tan, Z., Wang, Y.: Global well-posedness of an initial-boundary value problem for viscous non-resistive MHD systems. { SIAM J. Math. Anal.} {\bf 50}, 1432-1470(2018)


\bibitem{WW} Wu, J., Wu, Y.: Global small solutions to the compressible 2D magnetohydrodynamic system without magnetic diffusion. {Adva. Math.} {\bf 310}, 759-888(2017)


\bibitem{XZ} Xu, L., Zhang, P.: Global small solutions to three-dimensional incompressible magnetohydrodynamical system. { SIAM J. Math. Anal.} {\bf 47}, 26-65(2015)

\bibitem{ZZ} Zhou, Y., Zhu, Y.: Global classical solutions of 2D MHD system with only magnetic diffusion on periodic domain. {J. Math. Phys.} {\bf59}, 081505(2018)

\end{thebibliography}
\end{document}